\RequirePackage{fix-cm}
\RequirePackage{amsmath}
\RequirePackage{accents}

\documentclass[12pt]{svjour3}       
\usepackage{graphicx}
\usepackage{amsmath}
\usepackage{color}
 \usepackage{xcolor}
 \usepackage{enumerate}
  \usepackage{mathrsfs} 
 \usepackage[document]{ragged2e}
 \usepackage{amssymb}
 \usepackage{mathtools}
 \usepackage{easybmat} 
 \usepackage{amsfonts}
 \usepackage{bm}
 \usepackage{textcomp}

 \usepackage{lipsum}
 \usepackage{afterpage} 
\usepackage{hyperref}

\usepackage{wasysym}
\usepackage[toc]{appendix} 
\usepackage{float}
\newcommand*{\dt}[1]{%
  \accentset{\mbox{\large\bfseries .}}{#1}}

  \let\dt\dot
  \newcommand*{\Scale}[2][4]{\scalebox{#1}{$#2$}} 
  \usepackage[utf8]{inputenc}
\usepackage[english]{babel}

\begin{document}

\title{A Continuum Description of Failure Waves
}


\author{Hamid Said         \and
        James Glimm} 
\institute{H. Said \at
               Department of Applied Mathematics and Statistics, \\ Stony Brook University, Stony Brook, NY 11794, USA \\
              Tel.: +1-484-844-7636\\
              Fax: N/A\\
              \email{hsaid@ams.sunysb.edu}           
           \and
           J. Glimm \at
              Department of Applied Mathematics and Statistics,\\ Stony Brook University, Stony Brook, NY 11794, USA
}
\date{}
\maketitle

\begin{abstract}
\justify
Shattering of a brittle material such as glass occurs dynamically through a propagating failure wave, which however, can not be assigned to any of the classical waves of the elasto-plastic theories of materials. Such failure waves have been a topic of research for decades. In this paper, we build a thermodynamically consistent theory based on the idea that a failure wave is analogous to a deflagration wave. Our theory admits, as special cases, the classical models of Feng and Clifton. Two fundamental thermodynamic functions (the free energy and the entropy production rate) form the basis of our theory. Such a two-function approach allows for the construction of a new variational principle and a new Lagrangian formulation that produce the equations of motion. Finally, a linearization of this theory is examined to gain insight into the coupling between the diffusive and elastic wave phenomena. 
\keywords{brittle fracture \and failure wave \and internal variable \and dissipative systems \and variational principle}
\end{abstract}

\section{Introduction}
\label{intro}

\begin{flushleft}
\justify
Fracture is one of the two main modes of material damage. Brittle fracture and comminuted fracture waves (also called failure waves) are an extreme version of this phenomena, arising in engineering applications such as modeling low-frequency earthquake \cite{sapozhnikov2007mechanistic}, armor damage prevention \cite{walley2013introduction}, determining the fracture strength of ceramics \cite{callister2011materials}, and developing blast technology in underground mining \cite{zhang2016rock}. Likewise, applications of brittle fracture arise in medicine: kidney stone treatment via shock wave lithography \cite{neuberg2006trigger}, quantifying comminution in bone fracture \cite{beardsley2005interfragmentary}, and modeling tooth structure and function \cite{lee2011fracture}. It is known that failure waves represent a dynamic wave mode of brittle elastic fracture that do not correspond to any of the classical elastic or inelastic waves of solid mechanics. Behind the wave front the material is in a comminuted or microfractured state \cite{summary2007,feng2000}. 
\end{flushleft}

\begin{flushleft}
\justify
Our main result is the derivation of failure wave modeling form a consistent themrodyanimic framework based on a diffusive process. Our equations include the classic work of Feng \cite{feng2000} and Clifton \cite{clifton1993analysis} as special cases. We identify Clifton's model as the zero dissipation limit endpoint in a one parameter family of Feng models. Our analysis suggests that anisotropy of diffusion (longitudinal vs transverse) is significant, which we argue is related to an observable quantity, the shard size of the comminuted rubble. Furthermore, we present a new data analysis corresponding to the width and rise time of the failure wave. 
\end{flushleft}

\begin{flushleft}
\justify
The dynamic response of glass to impact has been investigated by researchers for many years (see selected references \cite{feng2012experimental,summary2007,kanel20051,bourne96}). Experiments have conclusively demonstrated that once glass experiences impact (e.g. uni-axial compression) close to, but below the Hugoniot elastic limit, the sample initially deforms elastically at the shock front and subsequently fails dramatically through the formation of multiple microcracks into a comminuted (rubblized) zone. This delay in time is interpreted as the result of a self-propagating failure wave moving at subsonic speed behind the initial elastic shock. The material behind the failure wave suffers a total loss in tensile strength, and a large reduction in shear strength, while little to no change occurs in the longitudinal stress across the failure front. Researchers have used plate impact, rod impact and projectile penetration experiments to investigate the properties of failure fronts in various types of glass. 
\end{flushleft}

\begin{flushleft}
\justify
Various models have been developed to describe the failure wave: its formation, the speed at which the wave propagates, and the stability of the failure front. However, a consistent thermodynamic basis for the dynamics governing the failure wave is still lacking. Partom proposed a simple phenomological model, which relates rate of damage accumulation to the gradient of damage \cite{partom1998modeling}. Kanel \emph{et al.} also suggested a phenomological model where the failure wave is a moving control surface and the shear strength behind this surface is set to zero \cite{kanel1991failure}. A more fundamental model was proposed by Clifton \cite{clifton1993analysis} who described the failure wave as a propagating phase boundary within a system of conservation laws.  All these models do, in fact, provide multiple perspectives into the dynamics of the failure wave, and do agree, to varying degrees, with experimental observations. Here we model the propagation of the failure wave as analogous to a slow combustion process, a point of view presented by Feng \cite{feng2000} and later adopted by Chen \textit{et al.} \cite{chenfeng2003} (in collaboration with Feng). 
\end{flushleft}

\begin{flushleft}
\justify
In 2000, R. Feng \cite{feng2000} proposed a physically based diffusive field equation to describe the propagation of the failure wave. As he observed, since there is an increase in the lateral stress associated with the arrival of the failure wave, while the longitudinal stress remains unchanged, the propagation of a failure wave is not the result of momentum balance. Rather, it is the progressive percolation of microfissures that drives the propagation of the failure wave; hence, leading to a diffusive model. In fact, Feng observes that the propagation of the failure wave resembles that of a (subsonic) deflagration wave. He, therefore, writes a field equation of parabolic type for an unknown $V_d$:
\begin{equation} \label{eq:1.1} 
\dfrac{\partial V_d}{\partial t}=\nabla\cdot(\mathbf{D}(\textbf{x},t\textbf{)} \cdot \nabla V_d)+ \beta(\textbf{x},t)
\end{equation}
where $t \in \mathbb R^+$ and $\textbf{x} \in \mathbb R^3$ denote the time and space coordinates, and $\mathbf{D}$ and $\beta$ are the second order damage diffusivity tensor and evolution function, respectively. The variable $V_d$ is a damage related quantity: the increase in the specific volume of the material if it is damaged and then completely unloaded. Feng's model simulations were successful in reproducing the profile of the lateral stress measurements for different glass specimens. 
\end{flushleft}

\begin{flushleft}
\justify
The main critique of Feng's diffusive model is that the field equation \eqref{eq:1.1} is not derived from a consistent thermodynamic formulation. In this work we show that Feng's equation can be recovered in a theoretical derivation from thermodynamically consistent postulates. 
\end{flushleft}

\begin{flushleft}
\justify
The paper is organized as follows. In Section 2, we formulate the equations of motion and the constitutive theory based on two functions: the (Helmholtz) free energy and the rate of entropy production. This derivation is based on the principle of conservation of momentum, and the first and second law of thermodynamics. This two-function approach is related to but distinct from Biot's derivation of the equations of classical thermoelasticity \cite{biot1956thermo}. Biot's formulation \cite{biot1956thermo}, which in contrast to our theory, depends on a modified (Biot's) free energy, rather than the thermodynamically correct free energy. In our context, all the dynamics and thermodynamics of the system under study are contained in these two functions, and as such they specify the equations of motion, which are of hyperbolic-parabolic type. The evolution-diffusion equation are then linearized to obtain an anisotropic extension of Feng's model. In this context, we analyze the two distinct diffusion coefficients and compare our new theoretical predictions for the width and time rise of the failure waves to experiment. Lastly, we show that in the appropriate limit Feng's model reduces to Clifton's conservative model. In Section 3, we introduce a thermodyanimc potential as a sum of the Lagrangian and a dissipation function. The minimization of the sum leads to the equations of motion and constitutive relations. The equations of motion can thus be written as  Lagrange's equations with dissipation. In Section 4, we linearize the equations of motion, thereby illustrating the two components of the dissipative Lagrangian formalism, and offering insight into the nature of the interaction of the reversible effects upon the irreversible process in the system. Finally, Section 5 offers concluding remarks.
\end{flushleft}

\section{The equations of motion}
\label{sec:1}

\begin{flushleft}
\justify
We formulate the governing equations in Lagrangian coordinates as is customary in solid mechanics. See \cite{gurtinthermo,marsdenelasticity} for more details. Let the body $\mathcal{B}_0  \subset \mathbb R^n$ denote the undeformed reference configuration. A point $\mathbf{X} \in \mathcal{B}_0$ is called a  Lagrangian or material point. We represent the coordinates on $\mathcal{B}_0$ by $\{ X_i\}$, $i=1, 2, ..., n$. Furthermore, by defining a time-dependent diffeomorphism $\phi: \mathcal{B}_0 \times \mathbb R^+ \rightarrow \mathcal{B}_t \subset \mathbb R^n$ to be a motion of $\mathcal{B}_0$, we can represent an  Eulerian or spatial point occupied by	 $\textbf{X}$ at time $t$ by $\textbf{x}=\phi(\textbf{X},t)$, and so we have $\textbf{X}=\phi(\textbf{X},0)$.
\end{flushleft}

\begin{flushleft}
\justify
Define the displacement vector $\textbf{U} = \textbf{U}(\textbf{X},t)$ in Lagrangian coordinates by $$ \textbf{U}(\textbf{X},t) = \textbf{x} - \textbf{X}, \hspace{6mm} U_i(\textbf{X},t) =  x_i - X_i,\hspace{5mm}  i=1, 2 ,...,n$$
\end{flushleft} 
where $\textbf{x} = \phi(\textbf{X},t)$. 
\begin{flushleft}
\justify
Moreover, define the deformation gradient $\textbf{F}= \frac{\partial \textbf{x}}{\partial \textbf{X}}$, and so
\begin{equation} \label{eq:2.1}
\textbf{F} = \textbf{I} + \nabla\textbf{U}, \hspace{6mm} F_{ij} = \delta_{ij} + U_{i,j}, \hspace{5mm}  i,j=1, 2 ,...,n
\end{equation}

\end{flushleft}

\subsection{Constitutive theory}
\label{sec:2}

\begin{flushleft}
\justify
In this section we develop the thermodynamics of elasticity coupled to diffusive internal variables based on the work of Biot \cite{biot1956thermo}, Maugin, \cite{maugin1990infernal} and Ziegler \cite{ziegler2012introduction}. The derivation is based on the determining two fundamental thermodynamic functions. Such a two function approach was first proposed by Biot to obtain the equations of thermoelasticity \cite{biot1956thermo}. Biot observed that by minimizing the sum of two quantities -- a modified free energy (sometimes referred to as Biot's potential) and the dissipation function -- he was able to produce the classical equations of thermoelasticity. However, the relation of Biot's thermoelastic potential to fundamental thermodynamic principles is not clear. 
\end{flushleft}

\begin{flushleft}
\justify
For isothermal processes such as the one we are studying, the Helmholtz free energy is the appropriate fundamental thermodynamic potential representing the thermodynamic state of the system \cite{lebon2008understanding}, and minimizing the free energy produces the balance of forces for isothermal processes \cite{solidmechanics}. Moreover, in the case of failure waves, there is entropy production associated with the irreversible process of the growth and propagation of microcracks into the material \cite{maugin1990infernal}. In fact, as we will show, the rate of entropy production is proportional to the dissipation function. Therefore, the Helmholtz free energy and the dissipation function constitute the basis of our two-function approach. A more general treatment of this approach can be found in \cite{ziegler2012introduction}. We rely on Internal Variable Theory and classical theory of irreversible thermodynamics to construct our two functions. 
\end{flushleft}

\begin{flushleft}
\justify
Internal Variable Theory (IVT) has many applications including the description of elasto-plastic fracture and modeling viscoelastic behavior \cite{ziegler2012introduction,murakami2012continuum}. The dialated volume $V_d$ in Feng's model can be interpreted as an internal variable. We adopt this approach and introduce an internal variable $\mathbf{\Gamma}$. We assert that the propagation of the failure wave is represented by an evolution-diffusion equation for the internal variable $\mathbf{\Gamma}$. The tensorial nature of $\mathbf{\Gamma}$ and its interpretation as it pertains to the failure wave phenomenon will be discussed below. In keeping with our analogy to thermoelasticity and combustion theory, the diffusion of $\mathbf{\Gamma}$ can be viewed in some respects as analogous to the diffusion of heat in a body as briefly mentioned by Chen \emph{et al.} \cite{chenfeng2003}.  
\end{flushleft}

\begin{flushleft}
\justify
We first begin by considering the specific Helmholtz free energy denoted $\Psi$, and assume the system to be at constant temperature $\Theta_0$. 
\end{flushleft}

\begin{flushleft}
\justify
In the theory of thermoelasticity, the Helmholtz free energy $\Psi$ is a function of the deformation gradient $\textbf{F}$ and temperature ${\Theta}$ such that \cite{gurtinthermo} $$ \textbf{S} = \rho_0 \frac{\partial \Psi}{\partial \textbf{F}}   \hspace{4.5mm} \textrm{and} \hspace{4.5mm} s = - \frac{\partial \Psi}{\partial \Theta} \ ,$$
where $\textbf{S}$ is the first Piola-Kirchhoff stress tensor, $s$ is the entropy, and $\rho_0$ is the mass density of the undeformed material.\\
In our context we assume that $\Psi$ depends on $\textbf{F}$ and $\mathbf{\Gamma}$ with 
\begin{equation} \label{eq:2.2}
\textbf{S} = \rho_0 \frac{\partial \Psi}{\partial \textbf{F}} \ .
\end{equation}
Note that $$\begin{aligned}
  S_{ij} = \rho_0\dfrac{\partial \Psi}{\partial F_{ij}} &= \rho_0\dfrac{\partial \Psi}{\partial U_{k,l}} \dfrac{\partial U_{k,l}}{\partial F_{ij}}  \\
 &=\rho_0\dfrac{\partial \Psi}{\partial U_{k,l}} \dfrac{\partial F_{kl}}{\partial F_{ij}} \\
  &= \rho_0\dfrac{\partial \Psi}{\partial U_{k,l}} \delta_{ki} \delta_{lj}  \\
  &= \rho_0\dfrac{\partial \Psi}{\partial U_{i,j}}  = \rho_0\dfrac{\partial \Psi}{\partial (\nabla\textbf{U})} \ . \\
\end{aligned} $$
\end{flushleft}

\begin{flushleft}
\justify
Before examining the associated variable to $\mathbf{\Gamma}$, we make the following two key assumptions: (i) the Helmholtz free energy $\Psi$ not only depends on $\textbf{F}$ (or equivalently on $\nabla\textbf{U}$) and the internal variable $\mathbf{\Gamma}$, but also on the gradient of the internal variable $\nabla \mathbf{\Gamma}$ as done by \cite{maugin1990infernal} (ii) function $\Psi$ is convex is in $\nabla \mathbf{\Gamma}$ and $\mathbf{\Gamma}$. We motivate these assumptions below. For now we can write symbolically 
\begin{equation} \label{eq:2.3}
\Psi = \Psi(\nabla\textbf{U}, \mathbf{\Gamma},\nabla \mathbf{\Gamma}) \ .
\end{equation}
Therefore, the time derivative of $\Psi$ is written as
\begin{equation} \label{eq:2.4}
\rho_0 \dot{\Psi} = \textbf{S} : \vec{\dot{\overline{\nabla \textbf{U}}}} - \mathbf{Z}  \vec{\dot{\Gamma}} - \nabla\cdot (\textbf{B}  \vec{\dt{{\Gamma}}})   \ ,
\end{equation}
where we have defined
\begin{equation} \label{eq:2.5}
\begin{aligned}
 \textbf{Z} = -\rho_0\frac{\delta \Psi}{\delta \mathbf{\Gamma}} = \textbf{A} -\nabla \cdot \textbf{B}, \hspace{6mm} \textbf{A} = -\rho_0\frac{\partial\Psi}{\partial  \mathbf{\Gamma}}, \hspace{6mm}  \textbf{B}=-\rho_0\frac{\partial\Psi}{\partial (\nabla \mathbf{\Gamma})}\\
\end{aligned}
\end{equation}
We consider $\textbf{Z}$ to be the associated variable to $\mathbf{\Gamma}$. The variable $\textbf{Z}$ is interpreted as an internal force corresponding to the internal variable $\mathbf{\Gamma}$ \cite{maugin1990infernal}.
\end{flushleft}

\begin{flushleft}
\justify
The second law of thermodynamics (or equivalently the Clausius-Duhem inequality) determines physically permissible processes. In Lagrangian coordinates it reads \cite{solidmechanics}
\begin{equation} \label{eq:2.6}
\frac{d}{dt}\int_{\Omega}\rho_0sd\Omega + \int_{\partial \Omega} \mathscr{\dot{S}} d\textbf{a} \geq 0 \ ,
\end{equation}
where $s$ is the entropy (per unit mass) of the system, $\Omega$ is the domain of the system, and $\mathscr{\dot{S}} $ is the entropy flow vector on the boundary. As the problem is isothermal the only relevant entropy is that associated with $\mathbf{\Gamma}$. Locally, \eqref{eq:2.6} simplifies to
\begin{equation} \label{eq:2.7}
\rho_0\dot{s}+ \nabla\cdot\mathscr{\dot{S}} \geq 0 \ .
\end{equation}
\end{flushleft}

\begin{flushleft}
\justify
By substituting $\Psi = \mathcal{E} -s\Theta_0$ ($\mathcal{E}$ being the internal energy of the system), and the balance of energy equation: $ \rho_0 \dot{\mathcal{E}} = \textbf{S}:\vec{\dot{F}} $ into the above equation, and taking into consideration equations \eqref{eq:2.5} we obtain
\begin{equation} \label{eq:2.8}
\textbf{Z} \vec{\dot{\Gamma}} +  \nabla\cdot (\textbf{B} \vec{\dot{\Gamma}} + \Theta_0 \dot{\mathscr{S}}) \geq 0 \ .
\end{equation}
By making the choice 
\begin{equation} \label{eq:2.9}
\dot{\mathscr{S}} = - \dfrac{\textbf{B} \vec{\dot{\Gamma}}}{\Theta_0} \ ,
\end{equation}
we arrive at the fundamental inequality 
\begin{equation}\label{eq:2.110} 
\textbf{Z} \vec{\dot{\Gamma}} \geq 0 \ .
\end{equation}
\end{flushleft}

\begin{flushleft}
\justify
In general diffusion can be expressed in terms of a flux $\mathbf{h}$ across the boundary. For instance in thermomechanics, $\mathbf{h}$ is given by Fourer's Law, which in turn gives rise to an entropy flux $\dot{\mathscr{S}}$ defined as $\dot{\mathscr{S}} =  \mathbf{h} / \Theta $. 
\end{flushleft}

\begin{flushleft}
\justify
In the context of our problem, there is no heat exchange at the boundary. So the diffusion of the internal variable is associated with a flux other than heat flux across the boundary. Equation \eqref{eq:2.9} gives an explicit expression for this flux: it is the entropy flux associated with internal variable $\mathbf{\Gamma}$. Therefore, by analogy to heat flow across the boundary, we obtain an explicit expression for the total rate of energy released from the surface of the material during the propagation of the failure wave:
\begin{equation}\label{seq} 
H =\int_{\partial \Omega} \textbf{h} \cdot {\textbf{n}}\hspace{1mm} da = \int_{\partial \Omega} \textbf{B} \vec{\dot{\Gamma}} \cdot \textbf{n} \hspace{1mm} da \ .
\end{equation}
where $\textbf{n}$ is the unit outer normal. Note it was by assuming $\Psi = \Psi ( \cdots, \nabla \mathbf{\Gamma})$ that we got $\textbf{B} \neq 0$, which gives a nonzero entropy flux $\dot{\mathscr{S}}$. 
\end{flushleft}

\begin{flushleft}
\justify
Had we only assumed $\Psi = \Psi (\nabla \textbf{U}, \mathbf{\Gamma})$, then we would expect $\dot{\mathscr{S}} = 0$ and an equation similar to \eqref{eq:2.110} would still hold (by replacing $\textbf{Z}$ with $\textbf{A}$) \cite{ziegler2012introduction,lebon2008understanding}. Thus, there still could be dissipation in the medium, but not diffusion of the internal variable. 
\end{flushleft}

\begin{flushleft}
\justify
Next we consider the dissipation function. It is known that the Helmholtz free energy determines the reversible processes occurring in the system. If in addition to reversible processes, the evolution of the system under consideration produces entropy (i.e. an accompanied irreversible process), then one needs to supplement the free energy with a scalar dissipation function $\mathcal{D}$. This approach has been used in thermomechanics by Ziegler \cite{ziegler2012introduction} and Biot \cite{biot1955variational,biot1956thermo}, and in a more general context by Maugin\renewcommand*{\thefootnote}{\arabic{footnote}}\footnote{In fact the dissipation function was first introduced by Relyigh \cite{rayleigh1945theory} in his study of dissipative processes in classical mechanics. He required such a function to be non-negative and quadratic in the velocities.} \cite{maugin1992thermomechanics}. 
\end{flushleft}

\begin{flushleft}
\justify
The constitutive law for the dissipation function $\mathcal{D}$ varies according to the system under study. In applications function $\mathcal{D}$ is positive and depends on the rate variables, as well as on the thermodynamic state variables \cite{ziegler2012introduction}. As a postulate, we exclude viscous related strain rate $\vec{\dot{\overline{\nabla U}}}$ from $\mathcal{D}$. Moreover, as the system is isothermal the only remaining rate variable is $\vec{\dot{\Gamma}}$, and so we write 
\begin{equation} \label{eq:2.11}
\mathcal{D} = \mathcal{D}(\vec{\dot{\Gamma}},w) \ ,
\end{equation}
where $w$ represents the state variables.
\end{flushleft}

\begin{flushleft}
\justify
Before investigating the functional form of $\mathcal{D}$ in the case of failure waves, we need to determine the forces associated with $\vec{\dot{\Gamma}}$. Following Ziegler \cite{ziegler2012introduction} and Maugin \cite{maugin1990infernal}, we assume that the function $\mathcal{D}$ is such that

\begin{equation} \label{eq:2.12}
\frac{\partial \mathcal{D}}{\partial \vec{\dot{\Gamma}}} + \frac{\delta \Psi}{\delta \mathbf{\Gamma}} = 0 \ ,
\end{equation}
or equivalently 
\begin{equation} \label{eq:2.13}
\rho_0\frac{\partial \mathcal{D}}{\partial \vec{\dot{\Gamma}} } = \textbf{Z} \ .
\end{equation}
In other words, the dissipative force associated with $\vec{\dot{\Gamma}}$ through $\mathcal{D}$ is identified as the opposite of the quasiconservative force associated with $\mathbf{\Gamma}$ through $\Psi$\footnote{According to Ziegler this is a statement that the net internal forces associated to the internal variables in an arbitrary process is equal to zero, which is a consequence of energy conservation.}.
\end{flushleft}

\begin{flushleft}
\justify
For failure waves, we assume the dependence of $\mathcal{D}$ on $w$ in \eqref{eq:2.11} vanishes, and furthermore, we assume a qaudratic  dependence of $\vec{\dot{\Gamma}}$ similar to the themroelastic theory of Biot \cite{biot1956thermo}
\begin{equation} \label{eq:2.14}
\rho_0\mathcal{D} = \rho_0 \mathcal{D}(\vec{\dot{\Gamma}}) =\dfrac{1}{2}\mathbf{\Lambda}\cdot( \vec{\dot{\Gamma}}\cdot \vec{\dot{\Gamma}})\geq 0 \ .
\end{equation}
For example if $\mathbf{\Gamma}$ is an $n$-vector, then the above equation in component form reads
\begin{equation} \label{eq:2.15}
\hspace{28mm} \rho_0\mathcal{D}(\vec{\dot{\Gamma}}) = \dfrac{1}{2}\Lambda_{ij} \dot{\Gamma}_i\dot{\Gamma}_j \geq 0\hspace{10mm}  i,j=1, 2 ,...,n \ .
\end{equation}
\end{flushleft}

\begin{flushleft}
\justify
The classical theory of irreversible thermodynamics relates the function $\mathcal{D}$ to the rate of production on internal entropy $\dot{s}^{(i)}$. Onsager's principle defines $\dot{s}^{(i)}$ in terms of fluxes and forces associated to the fluxes. In the absence of thermal effects the flux and force are identified with $\vec{\dot{\Gamma}}$ and $\mathbf{Z}$, respectively. So we write \cite{solidmechanics}
\begin{equation} \label{eq:2.16}
\rho_0 \Theta_0\dot{s}^{(i)} = \vec{\dot{\Gamma}} \cdot \mathbf{Z} \geq 0	 \ .
\end{equation}
If the fluxes and forces are related by a linear phenomenological relation 
\begin{equation} \label{eq:2.17}
\mathbf{Z} = \mathbf{\Lambda} \cdot \vec{\dot{\Gamma}}\ ,
\end{equation}
then,
\begin{equation} \label{eq:2.18}
\rho_0\dot{s}^{(i)} =\frac{1}{\Theta_0}  \vec{\dot{\Gamma}} \cdot \mathbf{Z}=\frac{1}{\Theta_0} \mathbf{\Lambda}\cdot \vec{\dot{\Gamma}} \cdot \vec{\dot{\Gamma}} = \dfrac{2 \rho_0}{\Theta_0}\mathcal{D} \ ,
\end{equation}
where by Onsager reciprocal relations we have $\mathbf{\Lambda} = \mathbf{\Lambda}^T$. Furthermore, we assume $\mathbf{\Lambda} $ is invertible. So not only does the dissipation produce entropy, the rate of entropy production is in fact proportional to the dissipation function. Hence, we can treat $\mathcal{D}$ as the source of entropy production. A more general treatment of dissipation functions can be found in \cite{ziegler2012introduction,gurtinthermo}.
\end{flushleft}

\begin{flushleft}
\justify
Expanding equation \eqref{eq:2.12} gives us
\begin{equation} \label{eq:2.19}
\frac{\partial \mathcal{D}}{\partial \vec{\dot{\Gamma}}}  + \frac{\partial \Psi}{\partial \mathbf{\Gamma}} -  \left( \frac{\partial \mathbf{P}}{\partial (\nabla \mathbf{U})} \cdot \nabla_{\textbf{X}}(\nabla \textbf{U}) + \frac{\partial \mathbf{P}}{\partial \mathbf{\Gamma}} \cdot \nabla_{\textbf{X}} \mathbf{\Gamma} + \frac{\partial \mathbf{P}}{\partial (\nabla \mathbf{\Gamma})} \cdot \nabla_{\textbf{X}}(\nabla \mathbf{\Gamma}) \right ) =0 \ ,
\end{equation}
where $\textbf{P} \doteq \rho_0 \dfrac{\partial \Psi}{\partial (\nabla \mathbf{\Gamma})} $ and symbol $\nabla_{\textbf{X}}$ means derivative with respect to $\textbf{X}$ while holding other variables constant. Moreover, we define the damage diffusion tensor $\mathbf{K}$ as
\begin{equation}\label{eq:2.20}
\mathbf{K} = \frac{\partial \mathbf{P}}{\partial (\nabla \mathbf{\Gamma})}=\rho_0 \dfrac{\partial^2 \Psi}{\partial (\nabla \mathbf{\Gamma}) \partial (\nabla \mathbf{\Gamma})} \ .
\end{equation}
\end{flushleft}

\begin{flushleft}
\justify
Tensor $\mathbf{K}$ is defined analogously to the classical elasticity tensor, and hence it is symmetric. Furthermore, since $\Psi$ is convex in $\nabla \mathbf{\Gamma}$ tensor $\mathbf{K}$ must be positive definite, thereby ensuring well-defined diffusion. 
\end{flushleft}

\begin{flushleft}
\justify
Before we formulate the hyperbolic-parabolic equations we discuss the order and nature of the variable $\mathbf{\Gamma}$. In general the order and physical meaning of an internal variable depends on the system under study. In the classical theory of material damage, a scalar variable $d \in \left[0,1\right]$ called a \textit{damage variable} is introduced such that $d=0$ corresponds to intact material, and $d=1$ corresponds to fully fractured material \cite{maugin1992thermomechanics}. However, in the context of our problem, the incident shock wave introduces lateral shear in the brittle material because of Poisson's ratio, which is released by the arrival of the failure wave as a drop in shear strength. It is reasonable then to assume that the difference in shear stress is driving the damage from a mechanical point of view. So it is natural to introduce a second order tensor to describe the diffusion of the damage $\mathbf{\Gamma}$. Variable $\mathbf{\Gamma}$ represents the strain due to the 3D anisotropic fracture and microcracks. The associated force (in this case a measure of shear stress) is given by the first equation of \eqref{eq:2.5}.
\end{flushleft}

\begin{flushleft}
\justify
For notational simplicity we rewrite the symmetric rank two indices $i, j$ of $\mathbf{\Gamma}$ as a single index $i =1, \ldots, \frac{n(n+1)}{2}$. This choice yields a rank four tensor $\textbf{K}$ and rank two tensor $\mathbf{\Lambda}$. 
\end{flushleft}

\subsection{Coupled hyperbolic-parabolic equations}
\label{sec:3}
\begin{flushleft}
\justify
We are now in position to write down the system of partial differential equations of hyperbolic-parabolic type describing the formation and propagation of failure waves
\begin{subequations} \label{eq:2.25}
\begin{align}
\rho_0\dfrac{\partial^2 \textbf{U}}{\partial t^2} - \nabla\cdot\textbf{S}\mid_{\mathbf{\Gamma = 0}} \hspace{1mm} &= \rho_0\mathbf{r} \label{eq:2.25a} \\
\frac{\partial \mathcal{D}}{\partial \vec{\dt{\Gamma}}} + \frac{\delta \Psi}{\delta \mathbf{\Gamma}} &= 0 \ , \label{eq:2.25b}
\end{align}
\end{subequations}
where $\textbf{r}(\textbf{X},t)$ is a (specific) body force. Equation \eqref{eq:2.25a} describes the motion of the initial nonlinear shock, hence we have a vanishing $\mathbf{\Gamma}$. 
\end{flushleft}

\begin{flushleft}
\justify
By using the constitutive equations for functions $\Psi$ and $\mathcal{D}$, particularly equations \eqref{eq:2.2}, \eqref{eq:2.15}, and \eqref{eq:2.20}, the coupled system\footnote{Strictly speaking there is one-way coupling only: the parabolic part is coupled with the hyperbolic equation since the hyperbolic part is independent of variable $\mathbf{\Gamma}$.} \eqref{eq:2.25} can be rewritten as

\begin{subequations} \label{eq:2.26}
\begin{eqnarray}
\rho_0\dfrac{\partial^2 U_i}{\partial t^2} - C_{ijkl} \dfrac{\partial^2U_k}{\partial X_j\partial X_l} & = & \rho_0 r_i \label{eq:2.26a} \\
\Lambda_{ij}\frac{\partial \Gamma_j}{\partial t} - \dfrac{\partial }{\partial X_j}\left( K_{ijkl} \dfrac{\partial \Gamma_k}{\partial X_l}\right) & = & f_i \left(U_{m,n}, \Gamma_m,  \Gamma_{m,n} \right) \ ,\label{eq:2.26b}
\end{eqnarray}
\end{subequations}
where $\textbf{C}\mid_{\mathbf{\Gamma = 0}} = \textbf{C}\left(\nabla \textbf{U},\mathbf{\Gamma} = 0\right)=\textbf{C}\left(\nabla \textbf{U}\right)$ is the elasticity tensor, $\textbf{f}$ is a nonlinear function in its arguments, and tensor $\mathbf{K}$ depends on the state vairables $\nabla \textbf{U}$, $\mathbf{\Gamma}$, $\nabla \mathbf{\Gamma}$. In component form the elasticity tensor is given to be
$$C_{ijkl} =\dfrac{\partial S_{ij}}{\partial U_{k,l}} \mid_{\mathbf{\Gamma = 0}}  =\rho_0 \dfrac{\partial^2\Psi}{\partial U_{i,j} \partial U_{k,l}} \mid_{\mathbf{\Gamma = 0}}  \ .$$
Hence, we have obtained a nonlinear coupled hyperbolic-parabolic system of PDEs. It should be realized that diffusion will be inactive (i.e. each term in equation \eqref{eq:2.26b} is identically zero) if the stress due to initial impact does not overcome a certain threshold $\sigma_0$. 
\end{flushleft}

\begin{flushleft}
\justify
System \eqref{eq:2.26} is constrained by the entropy inequality 	\eqref{eq:2.110} $$ \textbf{Z} \vec{\dt{\Gamma}} \geq 0 \ ,$$
which only allows for physically admissible solutions.   
\end{flushleft}

\begin{flushleft}
\justify
The initial and boundary conditions that must accompany the system of equations \eqref{eq:2.26} are given to be 
\begin{equation} \label{eq:2.27}
\textbf{U}(\textbf{X},0)=\textbf{U}^0(\textbf{X}), \hspace{2mm} \textbf{U}_t(\textbf{X},0)=\textbf{U}^1(\textbf{X}), \hspace{2mm} \mathbf{\Gamma}(\textbf{X},0)=\mathbf{\Gamma}^0(\textbf{X}), \hspace{5mm} \textbf{X} \in \mathcal{B}_0 \ ,
\end{equation}
\begin{equation} \label{eq:2.28}
\textbf{U}\mid_{\textbf{x} \in \partial \mathcal{B}_0}=\textbf{h}^1, \hspace{2mm} \mathbf{\Gamma}\mid_{\textbf{x} \in \partial \mathcal{B}_0}=\textbf{h}^2, \hspace{5mm} \forall t >0  \ ,
\end{equation}
where the functions  $\textbf{U}^0(\textbf{X}), \textbf{U}^1(\textbf{X}), \mathbf{\Gamma}^0(\textbf{X}), \textbf{h}^1, \textbf{h}^2 $ are all given. The dynamics \eqref{eq:2.26} of the problem is completely determined by the two functions: $\Psi$ and $\mathcal{D}$.
\end{flushleft}

\begin{flushleft}
\justify
In $\mathbb{R}^3$, the internal variable $\mathbf{\Gamma}$ has 6 independent components, and the rank 4 symmetric tensor $\mathbf{K}$ has 36 independent components and may be written as a $ 6 \times 6$ matrix. We obtain a simple non-isotropic model for diffusion if we assume that the material under study is isotropic, and consider a 1D shock wave in the longitudinal direction $X$; in such case the two lateral directions mirror one another. In addition, if we take a scalar variable $\Gamma$ to represent the damage, as done by Feng, then matrix $\textbf{K}$ reduces to a $ 3 \times 3$ matrix and the matrix $\mathbf{\Lambda}$ reduces to a scalar, which we denote $\lambda$. Finally, we assume that matrix $\textbf{K}$ is diagonalized in the longitudinal direction $X$ and the two transverse components. 
\end{flushleft}

\subsection{Feng's model}
\label{sec:4}

\begin{flushleft}
\justify
By the considerations given in the preceding paragraph, the unknown variables are written as $U = U(X,t)$ and $\Gamma =\Gamma (\textbf{X},t)$. We obtain Feng's diffusive field equation as follows. We retain the definition of $\mathcal{D}$ and expand \eqref{eq:2.3} to second order around the natural undamaged state (i.e. $ \nabla \textbf{U} = \mathbf{\Gamma} = \nabla \mathbf{\Gamma} = 0$). In this state the forces associated with state variables (i.e. equations \eqref{eq:2.2} and \eqref{eq:2.5}) vanish, and we assume that no coupling exists between the state variables. Thus the free energy and dissipation function are written as
\begin{equation} \label{eq:2.29} 
\rho_0\Psi_F = \dfrac{c_1}{2}U_{X}^2  + \frac{c_2}{2} \Gamma^2 + \frac{1}{2}K_{ij}  \nabla_i \Gamma \nabla_j \Gamma + \dfrac{c_3}{4}U_{X}^4 \hspace{4.5mm} \textrm{and} \hspace{5.5mm}  \rho_0\mathcal{D}_F = \dfrac{1}{2}\lambda \dt{\Gamma}^2 \ ,
\end{equation}
where constant $c_1$ is the 1D linear elastic coefficient (e.g. Young's modulus) and constant $c_3$ specifies the nonlinear elastic properties. An expression for scalar $c_2 = c_2(\textbf{X},t) \geq 0 $ is given by Feng \cite{feng2000} depending on the principle eigenvalue of $\mathbf{K}$ (see also \cite{chenfeng2003}). Furthermore, we allow $\textbf{K}$ to depend on variables $\textbf{X}$ and $t$.
\end{flushleft}

\begin{flushleft}
With these specifications equation \eqref{eq:2.26b} gives 

\begin{equation} \label{eq:2.37}
\lambda \frac{\partial \Gamma}{\partial t} - \dfrac{\partial }{\partial X_j}\left( K_{ij} \dfrac{\partial \Gamma}{\partial X_i}\right) =  - c_2\Gamma \ .
\end{equation} 
For a nonzero $\lambda$ \eqref{eq:2.37} can be written as

\begin{equation} \label{eq:2.30}
\frac{\partial \Gamma}{\partial t} - \dfrac{\partial }{\partial X_j}\left( D_{ij} \dfrac{\partial \Gamma}{\partial X_i}\right) =  c_4\Gamma \ ,
\end{equation}
where $c_4 = -\frac{c_2}{\lambda}$ and $D_{ij} = \frac{1}{\lambda} K_{ij}$:

\begin{equation}  \label{eq:2.31}
\textbf{D} = \begin{bmatrix}
d_1 & 0 \\ 
0 & d_2 \\
\end{bmatrix}
= \begin{bmatrix}
\frac{k_1}{\lambda} & 0 \\ 
0 & \frac{k_2}{\lambda}  \\
\end{bmatrix} \ ,
\end{equation}
where the coefficients $d_2$ and $k_2$ are multiples of the of the identity matrix $\textbf{I}_{2 \times 2}$.
\end{flushleft}

\begin{flushleft}
\justify
According to Feng $d_1$ determines the speed of the failure wave. We claim that $d_2$ determines the shard size in the damaged material. First, however, we discuss the evolution of $d_1$ and $d_2$ during the failure process. 
\end{flushleft}

\begin{flushleft}
\justify
In 2D projectile impact experiments in which a failure wave propagates as a spherical wave outwards, it is known \cite{wei2014mass} that after the abrupt cessation of the failure front isolated radial cracks start to form moving away from the impact point. We speculate that a similar phenomena for planar shock impact should occur parametrically as the shock impact is varied, with a transition from a failure wave to isolated cracks occurring at a critical shock strength.  To explain such behavior, we must have that during the propagation of the failure wave, the speed of diffusion associated with $d_2$, denoted $v_l$ (defined analoguosly to \eqref{eq:2.32}), is larger than the crack trip velocity $v_0$. As the failure wave propagates away from the impact point speed $v_l$ decreases until $v_l = v_0$ at which point the failure process stops, and the cracks ''escape'' the percolation in the lateral direction and as a result isolated radial cracks begin. 
\end{flushleft}

\begin{flushleft}
\justify
Chen \emph{et al.} \cite{chenfeng2003} observe that in order to obtain a planar failure wave in a uniaxially compressed material: "percolation
of microdamage [is] much faster in all the lateral (transverse) directions than in the longitudinal direction ... ". Hence, in terms of the diffusion coefficients we must have $ d_2 > d_1$ before the cessation of the failure wave. Chen \emph{et al.}, however, model only the isotropic case by assuming $d_1= d_2$, and they express $d_1$ in terms of a measure of shear stress. The point here is that according to their model the speed of the failure wave is governed only by $d_1$. 
\end{flushleft}

\begin{flushleft}
\justify
In diffusion-driven processes such as combustion, the speed of propagation is determined by the reaction rate --or equivalently reaction time-- and the diffusion coefficient \cite{clavin2016combustion}:
\begin{equation} \label{eq:2.32}
v_f \simeq \sqrt{d_1 \dot{\omega}} \hspace*{7mm} \textrm{or equivalently} \hspace*{7mm} v_f \simeq \sqrt{\frac{d_1}{\tau}}   \ ,
\end{equation}
where $\dot{\omega}$ and $\tau$ are the reaction rate and reaction time, respectively. By identifying $v_f$ with the speed of failure wave obtained by experiments, we can calculate the reaction time. Moreover, once $v_f$ is specified we are able to specify the width of the wave, $\delta_1$, through 
\begin{equation} \label{eq:2.33}
\delta_1 \doteq \dfrac{d_1}{v_f}  \ .
\end{equation}
 \end{flushleft}
 
 \begin{flushleft}
\justify
As a new analysis of existing experiments, we compare our theoretical predictions for the values of $\tau$ and $\delta_1$ for K8 glass and soda lime glass with the experiments examined by Feng. Experimental data are taken from plots of \cite{feng2000}. We have interpreted $\tau$ as the time rise in the lateral stress to find the experimental value of $\delta_1$.
\end{flushleft}

\begin{table}[H]
\centering
\caption{Comparison of model prediction with experimental data for K8 glass with $v_f = 3.32$ km/s, $d_1 = 6.6$ m$^2$/s. Gauge positioned at 4.5 mm from the loading surface.}
\begin{tabular}{c | c | c}
&
\hspace*{3mm} \footnotesize $\tau$ (s) \hspace*{3mm}&
\hspace*{3mm} \footnotesize $\delta_1$ (m) \hspace*{3mm} \\
\hline
&&\\
\footnotesize Experimental value & \footnotesize $0.75 \times 10^{-6}$ & \footnotesize $2.49 \times 10^{-3}$ \\
&&\\
\hline
&&\\
\footnotesize Theoretical prediction & \footnotesize $0.6 \times 10^{-6}$ & \footnotesize$2 \times 10^{-3}$ \\ 
&&\\
\end{tabular}
\end{table}

\begin{table}[H]
\caption{ Comparison of model prediction with experimental data for soda lime glass with $v_f = 3.09$ km/s, $d_1 = 7.4$ m$^2$/s. Gauge positioned at 3.3 mm from the loading surface.}
\centering
\begin{tabular}{c | c | c}
&
\hspace*{3mm} \footnotesize$\tau$ (s) \hspace*{3mm}&
\hspace*{3mm} \footnotesize$\delta_1$ (m) \hspace*{3mm} \\
\hline
&&\\
\footnotesize Experimental value & \footnotesize$0.9 \times 10^{-6}$ & \footnotesize$2.8 \times 10^{-3}$ \\
&&\\
\hline
&&\\
\footnotesize Theoretical prediction & \footnotesize$0.8 \times 10^{-6}$ & \footnotesize$2.4 \times 10^{-3}$ \\ 
&&\\
\end{tabular}
\end{table}  

\begin{flushleft}
\justify
Finally, we assert that the second eigenvalue $d_2$ of tensor $D_{ij}$ determines the fineness of fracture in the failed brittle material. If we define a measure of length $\delta_2$ in the following way $$\delta_2 \doteq \dfrac{d_2}{v_f}   \ , $$
then scalar $\delta_2$ gives an estimate for the fineness of fracture, which can be measured in principle via impact experiments.  
\end{flushleft}

\subsection{Clifton's model}
\label{sec:5}
\begin{flushleft}
\justify
The theory of deflagration can be formulated as a strictly hyperbolic system as found in \cite{courant1999supersonic,chorin1990mathematical}. However, if diffusion and finite reaction rates are considered, then combustion is obtained in the sense given by Zeldovich: a coupled hyperbolic-parabolic system of PDEs \cite{clavin2016combustion,bebernes2013mathematical}. In the latter formulation of the theory, the flame propagates with finite width and the flame speed depends on the diffusion coefficient and the reaction rate. In the absence of heat transfer and assuming an infinite reaction rate leading to a sharp flame front, Zeldovich's theory reduces to the hyperbolic system. This same sharp front limit can be achieved with zero viscous and thermal dissipations \cite{maugin1992thermomechanics,courant1999supersonic}. We show that in this zero dissipation limit for the propagation of the failure wave, Feng's model reduces to Clifton's conservative model. 
\end{flushleft}

\begin{flushleft}
\justify
Clifton \cite{clifton1993analysis} models the failure wave as a propagating phase boundary, which he calls a transformation shock. According to Clifton, the phase change occurs across a sharp front. Across the phase boundary, he imposes Rankine-Hugoniot conditions for the usual isentropic conservation laws in one space dimension, and derives an expression for the speed of the front. The task that lies before us, therefore, is to find the conditions under which the theory in Section 2.2 reduces to the conservation laws proposed by Clifton.
\end{flushleft}

\begin{flushleft}
\justify
Dissipation effects due to entropy production vanish for conservative systems, hence by taking $\mathbf{\Lambda} \longrightarrow 0$, we obtain $\mathcal{D} = 0$ in equation \eqref{eq:2.15} --assuming $\vec{\dot{\Gamma}}$ is bounded on $\mathcal{B}_0$. Moreover, since the entropy flux across the boundary $\dt{\mathscr{S}}$ must also vanish, we conclude $$ \mathbf{B} = - \rho_0 \dfrac{\partial \Psi}{\partial \nabla \mathbf{ \Gamma}} = 0 \ ,$$ which implies$$\mathbf{K} = 0 \ ,$$  and 
\begin{equation} \label{eq:2.34}
\frac{\partial s}{\partial t} = \dfrac{\partial s^{(i)}}{\partial t} = 0 \ .
\end{equation} 
Furthermore, by assuming the nonhomogenous term $f$ to be identically zero, equation \eqref{eq:2.26b} trivially vanishes. We claim that equations \eqref{eq:2.26a} and \eqref{eq:2.34}, produce the 1D conservation laws considered by Clifton.
\end{flushleft}

\begin{flushleft}
\justify
It is sufficient to show that \eqref{eq:2.34} is equivalent to the law of conservation of energy \cite{dafermos2005hyperbolic}
\begin{equation} \label{eq:2.35}
\dfrac{\partial}{\partial t} \rho_0\left( \mathcal{E} + \frac{1}{2}v^2 \right) - \frac{\partial}{\partial X} \left( vS \right) = 0 \ ,
\end{equation}
where we have defined $\mathcal{E} $ as the internal energy, $S$ is the scalar valued Piola-Kirchhoff stress, and $v = {\partial U}/{\partial t}$.
\end{flushleft}
\begin{flushleft}
\justify
Gibbs equation with an internal variable reads \cite{lebon2008understanding}
\begin{equation*}
d\mathcal{E} = \Theta_0 ds + \dfrac{1}{\rho_0}	SdU_X - Ad\Gamma \ ,
\end{equation*}
with $A = - \rho_0 \partial \Psi / \partial \Gamma$ in \eqref{eq:2.5}. However, by \eqref{eq:2.13} we have $A=0$. Thus, Gibbs equation reduces to
\begin{equation} \label{eq:2.36}
d\mathcal{E} = \Theta_0 ds + \dfrac{1}{\rho_0}SdU_X \ .
\end{equation}
Assuming $\mathcal{E} = \mathcal{E}(s,U_X)$ combined with \eqref{eq:2.36} we obtain $$ \dfrac{\partial \mathcal{E}}{\partial t} = \Theta_0 \frac{\partial s}{\partial t} + \dfrac{1}{\rho_0} S \frac{\partial U_X}{\partial t} = \Theta_0 \ . \frac{\partial s}{\partial t} +\dfrac{1}{\rho_0} S \frac{\partial v}{\partial X} \ .$$
Therefore $$\dfrac{\partial}{\partial t} {\rho_0}\left( \mathcal{E} + \frac{1}{2}v^2 \right) = {\rho_0}\dfrac{\partial  \mathcal{E} }{\partial t}  + {\rho_0} v\frac{\partial v}{\partial t} = {\rho_0}\Theta_0 \frac{\partial s}{\partial t} + S \frac{\partial v}{\partial X} + v \dfrac{\partial S }{\partial X} \ , $$
where we have used the balance of momentum equation $\rho_0\partial v /\partial t = \partial S / \partial X$. By substituting for \eqref{eq:2.34}, we obtain the conservation of energy equation \eqref{eq:2.35}.
\end{flushleft}

\begin{flushleft}
\justify
We have established the classical conservation laws
\begin{equation} \label{eq:2.38}
\begin{aligned}
& \hspace*{3mm} \dfrac{\partial U_X}{\partial t} - \dfrac{\partial v}{\partial X} = 0 \\
&\rho_0\dfrac{\partial v}{\partial t} -  \dfrac{\partial S	}{\partial X} = 0 \\
&\rho_0 \dfrac{\partial}{\partial t} \left( \mathcal{E} + \frac{1}{2}v^2 \right) - \frac{\partial}{\partial X} \left( vS \right) = 0  \ .
\end{aligned}
\end{equation}
\end{flushleft}

\begin{flushleft}
\justify
In one space dimension, if we choose $\Gamma = V_d$ and write Feng's equation in the form \eqref{eq:2.37}, then the limit $\lambda \longrightarrow 0$ implies $\mathcal{D}_F = 0$ and equation \eqref{eq:2.34} is satisfied. Moreover, we assume that the nonhomogeneous term is zero (i.e. $c_2 = 0$) as postulated earlier; then by the above argument, we obtain system \eqref{eq:2.38} and, hence, Feng's model reduces to the Clifton model.
\end{flushleft}

\begin{flushleft}
\justify
The field equations \eqref{eq:2.38} satisfy the Rankine-Hugoniot jump relations, which is the starting point considered by Clifton for modeling failure waves. We have, therefore, demonstrated that Clifton's model can be recovered from Feng's model in the limit $\mathbf{\Lambda} \longrightarrow 0$. In this limit the entropy is constant and, therefore, entropy production associated with the irreversible process of the growth of microcracks is zero over time, making Clifton's model an idealization of the process of the propagation of failure waves. 
\end{flushleft}

\section{Variational principle and Lagrangian formalism}
\label{sec:6}

\begin{flushleft}
\justify
In this section we formulate a a novel variational principle which produce the equations \eqref{eq:2.26}. We show that minimizing the sum of the Lagrangian function and the dissipation function leads to the equations of motion. Furthermore, we demonstrate that system \eqref{eq:2.26} can be rewritten as Lagrnage's equations with dissipation similar to dissipative systems in classical mechanics. 
\end{flushleft}

\begin{flushleft}
\justify
We start with a simplified variational principle where all surface and boundary terms vanish. While this is only a special case, the main physical and mathematical insights still hold. Appendix A presents the complete variational principle that includes all boundary terms, and is formulated in terms of generalized coordinates. 
\end{flushleft}

\begin{flushleft}
\justify
For the body $\mathcal{B}_0 $, define the total (Helmholtz) free energy $\psi$, the total kinetic energy $\mathscr{K}$, and the total dissipative function $\mathscr{D}$:

\begin{equation*}
\psi = \int_{\mathcal{B}_0}\rho_0 \Psi dV, \hspace{2.5mm} \mathscr{K} =\Scale[1.2]{\frac{1}{2}}\int_{\mathcal{B}_0}\rho_0 \frac{\partial U_i}{\partial t}\frac{\partial U_i}{\partial t}  dV, \hspace{2.5mm} \mathscr{D}=\Scale[1.2]{\frac{1}{2}}\int_ {\mathcal{B}_0} \Lambda_{ij} \dt{\Gamma}_i\dt{\Gamma}_j dV \ .
\end{equation*}

\end{flushleft}

\begin{flushleft}
\justify
Furthermore, we define the Lagrangian $\mathscr{L}$

\begin{equation*} 
\mathscr{L}  =\psi -  \mathscr{K}  = \mathscr{L}(\nabla \textbf{U}, \vec{\dt{U}} ,\mathbf{\Gamma}, \nabla \mathbf{\Gamma}) = \int_{\mathcal{B}_0} L dV \ .
\end{equation*}

\end{flushleft}

\begin{flushleft}
\justify
With the above definitions, the variational principle for the equations of motion is
 \begin{equation} \label{eq:3.1}
\delta I = \delta \int^{t_2}_{t_1}(\mathscr{L} + \mathscr{D})dt = 0 \ ,
\end{equation}
where $t_1$ and $t_2$ represent fixed instants of time. The variation is taken with respect to the displacement $\textbf{U}$ and the internal variable $\mathbf{\Gamma}$. 
\end{flushleft}

\begin{flushleft}
\justify
We demonstrate that the variational principle \eqref{eq:3.1} produces the system of equations \eqref{eq:2.26}. 
\end{flushleft}

\begin{flushleft}
\justify
We define the variation of the dissipation function $\delta \mathscr{D}$ in a similar fashion to Biot:  
\begin{equation} \label{eq:3.2}
\delta \mathscr{D} = \int_{\mathcal{B}_0} \Lambda_{ij}\dt{\Gamma}_i\delta \Gamma_j dV \ .
\end{equation}
The justification of equation \eqref{eq:3.2} becomes evident when we introduce generalized coordinates, as done in Appendix A.  
\end{flushleft}

\begin{flushleft}
\justify
Thus,
 \begin{align*}
&\int^{t_2}_{t_1}(\delta \mathscr{D} + \delta \psi)dt - \Scale[1.2]{\frac{1}{2}}\int^{t_2}_{t_1} \delta \int_{\mathcal{B}_0}\rho_0 \frac{\partial U_i}{\partial t}\frac{\partial U_i}{\partial t}  dV dt = 0 &\\
&\int^{t_2}_{t_1}\lbrace\delta \mathscr{D} +\int_{\mathcal{B}_0}\rho_0 \left(\frac{\delta \Psi}{\delta U_i}\delta U_i + \frac{\delta \Psi}{\delta \Gamma_i}\delta \Gamma_i\right) dV \rbrace dt - \int_{t_!}^{t_2}\int_{\mathcal{B}_0}\rho_0 \frac{\partial U_i}{\partial t}\frac{\partial (\delta U_i)}{\partial t}  dV dt = 0& \\
&\int^{t_2}_{t_1}\lbrace\delta \mathscr{D} +\int_{\mathcal{B}_0}\rho_0 \left(-\nabla_j\cdot \left( \frac{\partial \Psi}{\partial U_{i,j}} \right) \delta U_i + \frac{\delta \Psi}{\delta \Gamma_i}\delta \Gamma_i\right) dV \rbrace dt \\  &  \qquad \qquad \qquad \qquad \qquad \qquad \qquad \qquad \qquad \hspace{1mm}+ \int_{t_!}^{t_2}\int_{\mathcal{B}_0}\rho_0 \frac{\partial^2 U_i}{\partial t^2} \delta U_i dV dt = 0 &\\
&\int^{t_2}_{t_1}\lbrace \int_{\mathcal{B}_0} \left( -\nabla\cdot\textbf{S} +\rho_0 \frac{\partial^2\textbf{U}}{\partial t^2} \right) \delta \textbf{U} dV + \int_{\mathcal{B}_0} \left( \mathbf{\Lambda} \cdot \vec{\dt{\Gamma}} + \rho_0\frac{\delta \Psi}{\delta \mathbf{\Gamma}} \right) \delta \mathbf{\Gamma}dV \rbrace dt = 0 & \ . 
\end{align*}Since the last equation holds for all variations $ \delta \textbf{U}$ and $\delta \mathbf{\Gamma}$, we are left with the desired system of equations \eqref{eq:2.26} with $\textbf{r} = 0$.
\end{flushleft}

\begin{flushleft}
\justify
Therefore, the associated Lagrange equations to variational principle \eqref{eq:3.1}, which is equivalent to the governing system of PDEs \eqref{eq:2.26}, is written as
\begin{equation} \label{eq:3.3}
 \frac{\delta \mathscr{L}}{\delta \eta_i} - \frac{d}{dt}\left( \frac{\partial \mathscr{L}}{\partial \dt{\eta}_i} \right)  + \frac{\partial \mathscr{D}}{\partial \dt{\eta}_i } = 0 \ ,
\end{equation}
where we have defined the vector variable $\boldsymbol \eta = (\textbf{U},\mathbf{\Gamma})$. Notice that the form of equation \eqref{eq:3.3} is identical to the form of Lagrange's equations of motion in classical mechanics for dissipative systems \cite{goldstein1965classical}.
\end{flushleft}

\begin{flushleft}
\justify
There exists dissipative Hamiltonian and bracket formulations corresponding to the above dissipative Lagrangian system. The construction of said formalisms can be found in the thesis of the first author.

\end{flushleft}

\section{Linear theory}
\label{sec:6}
\begin{flushleft}
\justify
In this section we linearize the equations \eqref{eq:2.26} and note the interaction of the reversible effects upon the irreversible process of material failure.
\end{flushleft}

\begin{flushleft}
\justify
We begin by requiring that the free energy $\Psi$ be a quadratic form in $\nabla \textbf{U}$ and $\mathbf{\Gamma}$ and quadratic in $\nabla \mathbf{\Gamma}$ as follows
\begin{equation} \label{eq:4.1}
\rho_0\Psi = \frac{1}{2} \mathbf{a} \cdot \left( \nabla \mathbf{U} \cdot \nabla \mathbf{U} \right) + \textbf{b} \cdot \left( \nabla \mathbf{U} \cdot \mathbf{\Gamma} \right) + \frac{1}{2}\textbf{c} \cdot \left( \mathbf{\Gamma} \cdot \mathbf{\Gamma} \right)+ \frac{1}{2}\textbf{K} \cdot\left( \nabla \mathbf{\Gamma} \cdot \nabla\mathbf{\Gamma}\right)
\end{equation}
where $\textbf{a} $ is the constant (linear) elasticity tensor, $\textbf{b}$ is the constant coupling reversible effects to the irreversible ones, and $\textbf{c}$ is a dissipative constant. Moreover, we retain the definition of the quadratic dissipation function $\mathcal{D}$ as given in \eqref{eq:2.15}. We notice that by adding a nonlinear quartic term in $\nabla \mathbf{U}$ and setting $\textbf{b} =0$, we obtain the free energy for Feng's model in three dimensions (compare with \eqref{eq:2.29}). 

\end{flushleft}

\begin{flushleft}
\justify
With the above identifications, we substitute the Lagrangian function $$\mathscr{L} = \int_{\mathcal{B}_0}\rho_0 \left(  \Psi - \frac{1}{2} \frac{\partial U_i}{\partial t}\frac{\partial U_i}{\partial t} \right) dV \ ,$$ and dissipation function $\mathscr{D}$ into equations of motion \eqref{eq:3.3} to obtain
\begin{subequations} \label{eq:4.2}
\begin{align}
\rho_0\dfrac{\partial^2\mathbf{U}}{\partial t^2}  - \mathbf{a} \hspace{0.0mm} \bigtriangleup\hspace{-0.68mm}\mathbf{U} =0 \label{eq:4.2a}\\
\mathbf{\Lambda} \dfrac{\partial \mathbf{\Gamma}}{\partial t} - \nabla \cdot \left( \textbf{K} \hspace{1mm}  \nabla  \mathbf{\Gamma}\right) + \textbf{b} \hspace{0.5mm} \nabla \textbf{U} + \textbf{c} \hspace{0.5mm}\mathbf{\Gamma} = 0 \label{eq:4.2b}
\end{align}
\end{subequations}
where we have evaluated $\rho_0\frac{\partial \Psi}{\partial (\nabla\textbf{U})}$ at $\mathbf{\Gamma} = 0$ in equation \eqref{eq:4.2a}. In the above linear system of PDEs we have obtained a linear diffusion equation coupled to a  linear wave equation. 
\end{flushleft}

\begin{flushleft}
\justify
We notice that for a nonzero $\textbf{b}$, equation \eqref{eq:4.2b} picks up the reversible contribution $\textbf{b} \hspace{0.5mm} \nabla \textbf{U}$, as observed. Hence, for $\textbf{b} =0$, equation \eqref{eq:4.2b} is purely dissipative with no additional elastic effects, which is precisely what Chen \textit{et al.} consider in their model \cite{chenfeng2003}: "[T]he failure wave behind the shock wave front will not cause additional strain at the macroscopic level under the plate impact conditions ...". 
\end{flushleft}

\begin{flushleft}
\justify
The constant tensor $\textbf{b}$, therefore, measures the \emph{coupling} of the irreversible effects to the reversible effects. One expects that in an ideal situation constant $\textbf{b}$ should be set to zero, but due to varying material properties and/or different experimental configurations, residual elastic effects in the dissipative process may appear as a result, hence a non-zero $\textbf{b}$. This observation explains the inconsistencies in the data observed by Feng as they pertain to the rise in the longitudinal strains behind the failure wave \cite{feng2000}.
\end{flushleft}

\section{Conclusions}
\label{sec:7}
\begin{flushleft}
\justify
We have formulated a theory describing failure waves in brittle elastic material at a continuum level based on a thermodynamically consistent theory. We have subsumed and extended in our analysis a prior model proposed by Feng for brittle fracture, and we have recovered Clifton's model in the dissipationless limit. Our analysis reveals the importance of the coefficient of lateral diffusion $d_2$ in determining both the shard size and the transition from failure waves to isolated cracks. In this work we have combined and modified several methodologies previously developed in the context of thermoelascticty, IVT, and irreversible thermodynamics. 
\end{flushleft}

\begin{flushleft}
\justify
In summary, we developed a constitutive theory and the equations of motion in the context of two fundamental thermodynamic functions: the (Helmholtz) free energy $\Psi$ and a dissipation function $\mathcal{D}$. Feng's model is recovered by specific choices of the potentials $\Psi = \Psi _F$ and $\mathcal{D} = \mathcal{D}_F$ defining the constitutive theory. Clifton's model, in turn, can be recovered from Feng's by taking the limit $\mathbf{\Lambda} \longrightarrow 0$. Our two-function approach gives rise to a variational principle and a Lagrangian formalism. Finally, we presented a linear theory to gain insight into the  interactions of the reversible and irreversible processes.\\ 
\end{flushleft}

\begin{acknowledgements}
\justify
The authors thank Roman Samulyak for his helpful comments. Hamid Said wishes to thank the generous support of Kuwait University through the graduate scholarship. James Glimm's research is supported in part by the Army Research Organization, award \#W911NF1310249.\\
\end{acknowledgements}

\clearpage

\nocite{*}
\bibliographystyle{abbrv}
\bibliography{unt.bbl}


%
%
\clearpage

\appendix
\section{Appendix}
\begin{flushleft}
\justify
In this appendix we derive Lagrange's equations \eqref{eq:3.3} in generalized coordinates from a variational principle that includes boundary terms. We adapt Biot's method \cite{biot1970variational} for formulating a variational principle for thermoelasticity to our problem.
\end{flushleft}

\begin{flushleft}
\justify
We begin by stating the complete variational principle:
\begin{equation} \label{eq:A.1}
\delta \int^{t_2}_{t_1}(\psi - \mathscr{K} + \mathscr{D})dt = \int^{t_2}_{t_1} R   dt \ ,
\end{equation}
where 
\begin{equation*}
R= \int_{\partial \mathcal{B}_0} \left(S_{ij}n_j \delta U_i - B_{ij} n_j \delta \Gamma_i \right) dA + \int_{\mathcal{B}_0} \rho_0 r_i \delta U_i dV \ .
\end{equation*}
\end{flushleft}

\begin{flushleft}
\justify
We start by expressing the unknown vectors $\textbf{U}$ and $\mathbf{\Gamma}$ in terms of $k$ generalized coordinates $\mathbf{q} = (q_1, ..., q_k)$. We write
\begin{equation} 
\mathbf{U} = \mathbf{U} ( q_1, ..., q_k, \textbf{X},t), \hspace{9mm} \mathbf{\Gamma} = \mathbf{\Gamma} ( q_1, ..., q_k, \textbf{X},t) \ .
\end{equation}
\end{flushleft}

The variations of variables $\textbf{U}$ and $\mathbf{\Gamma}$ are given in terms of variation $\delta \mathbf{q}$:
\begin{equation}
\delta{U}_i = \dfrac{\partial U_i}{\partial q_j} \delta q_j , \hspace{9mm} \delta  \Gamma_i = \dfrac{\partial \Gamma_i}{\partial q_j} \delta q_j  \ .
\end{equation}

First, we look at $\int_{\mathcal{B}_0}\Lambda_{ij}\dt{\Gamma}_i\delta \Gamma_j dV $:

\begin{align*}
\int_{\mathcal{B}_0}\Lambda_{ij}\dt{\Gamma}_i\delta \Gamma_j dV 
&= \int_{\mathcal{B}_0}\Lambda_{ij}\dt{\Gamma}_i \dfrac{\partial {\Gamma}_j}{\partial {q}_m} \delta q_m dV \\
&= \int_{\mathcal{B}_0}\Lambda_{ij}\dt{\Gamma}_i \dfrac{\partial \dt{\Gamma}_j}{\partial \dt{q}_m} \delta q_m dV \hspace{7mm} (\textrm{for} \hspace{1mm} \frac{\partial {\Gamma}_j}{\partial {q}_m} = \frac{\partial \dt{\Gamma}_j}{\partial \dt{q}_m}) \\
 &= \dfrac{\partial}{\partial \dt{q}_m} \int_{\mathcal{B}_0}\dfrac{1}{2}\Lambda_{ij}\dt{\Gamma}_i \dt{\Gamma}_j \delta q_m dV \\ \ .
\end{align*} 
So if we define
\begin{align}
 \mathscr{D}=\Scale[1.2]{\frac{1}{2}}\int_ {\mathcal{B}_0} \Lambda_{ij} \dt{\Gamma}_i\dt{\Gamma}_j dV \ ,
 \end{align}
we obtain
\begin{align}
\int_{\mathcal{B}_0}\Lambda_{ij}\dt{\Gamma}_i\delta \Gamma_j dV = \dfrac{\partial \mathscr{D}}{\partial \dt{q}_m}\delta q_m  \doteq \delta \mathscr{D}\ .
\end{align}

\begin{flushleft}
\justify
Next, we examine the variation of $\psi$. We assume that $\psi$ is independent of the generalized velocities $\vec{\dt{q}}$. Moreover, below we keep the boundary terms when integrating by parts

\begin{align*}
\delta \psi &= \rho_0\int_{\mathcal{B}_0} \left( \dfrac{\partial \Psi}{\partial U_{i,j}}\delta U_{i,j} + \dfrac{\partial \Psi}{\partial \Gamma_{i}}\delta \Gamma_{i} + \dfrac{\partial \Psi}{\partial \Gamma_{i,j}}\delta \Gamma_{i,j}  \right) dV\\
&= \rho_0\int_{\mathcal{B}_0} \left( - \nabla_j \cdot \left( \dfrac{\partial \Psi}{\partial U_{i,j}} \right) \delta U_i +\dfrac{\partial \Psi}{\partial \Gamma_{i}}\delta \Gamma_{i} - \nabla_j \cdot \left( \dfrac{\partial \Psi}{\partial \Gamma_{i,j}} \right) \delta \Gamma_i\right) dV \\
&\textcolor{white}{+} \qquad \qquad \qquad\qquad\qquad\qquad+\rho_0\int_{\partial \mathcal{B}_0} \left( \dfrac{\partial \Psi}{\partial U_{i,j}}n_j\delta U_i + \dfrac{\partial \Psi}{\partial \Gamma_{i,j}}n_j\delta \Gamma_i \right) dA \\
&= \int_{\mathcal{B}_0} \left( -S_{ij,j} \delta U_i + \rho_0\dfrac{\delta \Psi}{\delta \Gamma_i}\delta \Gamma_i \right) dV + \int_{\partial\mathcal{B}_0} \left(T_i\delta U_i - B_{ij}n_j\delta \Gamma_i \right) dA \\
&= \dfrac{\partial \psi}{\partial q_m}\delta q_m \ ,
\end{align*} 

with $T_i = S_{ij}n_j$ is the force per unit area applied at the boundary of $\mathcal{B}_0$, and $B_{ij} = -\rho_0\dfrac{\partial \Psi}{\partial \Gamma_{i,j}}$ is associated with the entropy flux (see Section 2.1). In addition, we define the generalized force $Q_i$ as $$ Q_i = \int_{\partial\mathcal{B}_0} \left(T_j \dfrac{\partial U_j}{\partial q_i}  - B_{lj}n_j\dfrac{\partial \Gamma_l}{\partial q_i}  \right) dA + \int_{\mathcal{B}_0}\rho_0 r_j \dfrac{\partial U_j}{\partial q_i} dV$$
\end{flushleft}
where $\textbf{r}(\textbf{X},t)$ is a body force (see Section 2.2), so that $$R =Q_i \delta q_i $$
\begin{flushleft}
\justify
Lastly. we look at the variation of the kinetic energy $\mathscr{K}$. In Section 3 we calculated $\delta \mathscr{K}$
\begin{equation*}
-\delta \mathscr{K} = \int_{\mathcal{B}_0} \rho_0 \ddot{U}_i \delta U_i dV = \int_{\mathcal{B}_0} \rho_0 \ddot{U}_i \dfrac{\partial U_i}{\partial q_m}\delta q_m dV
\end{equation*}
Since $\frac{\partial U_i}{\partial q_m} = \frac{\partial \dt{U}_i}{\partial \dt{q}_m}$, we obtain
\begin{equation*}
\ddot{U}_i \dfrac{\partial U_i}{\partial q_m} = \dfrac{d}{dt} \left( \dt{U}_i \frac{\partial \dt{U}_i}{\partial \dt{q}_m} \right) - \dt{U}_i \dfrac{\partial \dt{U}_i}{\partial q_m}
\end{equation*}
Thus we arrive at
\begin{equation*}
- \delta \mathscr{K} = \int_{\mathcal{B}_0} \rho_0 \ddot{U}_i \delta U_i dV = \left\lbrace \dfrac{d}{dt} \left( \dfrac{\partial \mathscr{K}}{\partial \dt{q}_m} \right)- \dfrac{\partial \mathscr{K}}{\partial q_m} \right \rbrace \delta q_m
\end{equation*}

\end{flushleft}

\begin{flushleft}
\justify
Finally, by substituting the above three variations into the variational principle \eqref{eq:A.1} we obtain
\begin{equation*}
\delta \int^{t_2}_{t_1}(\psi - \mathscr{K} + \mathscr{D})dt = \int^{t_2}_{t_1}\left( \frac{\partial \psi}{\partial q_i} + \frac{d}{dt} \left( \frac{\partial \mathscr{K}}{\partial \dt{q}_i} \right) - \frac{\partial \mathscr{K}}{\partial q_i}   + \frac{\partial \mathscr{D}}{\partial \dt{q_i}} \right) \delta q_i dt = \int_{t_1}^{t_2} Q_i \delta q_i dt
\end{equation*}
Since the time interval $[t_1,t_2]$ is arbitrary we are left with 
\begin{equation} 
\dfrac{d}{dt} \left( \dfrac{\partial \mathscr{K}}{\partial \dt{q}_i} \right)- \dfrac{\partial \mathscr{K}}{\partial q_i} + \dfrac{\partial \mathscr{D}}{\partial \dt{q}_i} + \dfrac{\partial \psi}{\partial q_i} = Q_i
\end{equation} 
\end{flushleft}

\begin{flushleft}
\justify
We can rewrite equation (A.6) in the form of Lagrange's equations with dissipation by first introducing the full potential energy of the system $$ \mathcal{V} = \int_{\mathcal{B}_0} \rho_0 \Psi dV + \int_{\mathcal{B}_0}\rho_0 r_j  U_j dV + \int_{\partial\mathcal{B}_0} \left(T_i U_i - B_{ij}n_j\Gamma_i\right) dA \ ,$$ 
and defining the Lagrangian $\mathscr{L} \doteq \mathcal{V} - \mathscr{K} $:
\begin{equation} 
\frac{\delta \mathscr{L}}{\delta q_i} -  \frac{d}{dt}\left( \frac{\partial \mathscr{L}}{\partial \dt{q}_i} \right) + \frac{\partial \mathscr{D}}{\partial \dt{q}_i } = 0
\end{equation}
This concludes the derivation.
\end{flushleft}

\begin{flushleft}
\justify
We note that by substituting the expression for the variations of $\mathscr{\psi}$, $\mathscr{K}$, and $\mathscr{D}$ into the variational principle \eqref{eq:A.1} we obtain the full nonlinear system \eqref{eq:2.26}. 
	
\end{flushleft}

\end{document}